\documentclass{amsart}

\newtheorem{theorem}{Theorem}
\newcommand{\bt}{\begin{theorem}}
\newcommand{\et}{\end{theorem}}

\newtheorem*{theoremNN}{Theorem}
\newcommand{\btNN}{\begin{theoremNN}}
\newcommand{\etNN}{\end{theoremNN}}

\newtheorem{lemma}{Lemma}
\newcommand{\bl}{\begin{lemma}}
\newcommand{\el}{\end{lemma}}

\newtheorem{corollary}{Corollary}
\newcommand{\bc}{\begin{corollary}}
\newcommand{\ec}{\end{corollary}}

\newtheorem{definition}{Definition}
\newcommand{\bdf}{\begin{definition}}
\newcommand{\edf}{\end{definition}}

\newtheorem{conjecture}{Conjecture}
\newcommand{\bconj}{\begin{conjecture}}
\newcommand{\econj}{\end{conjecture}}

\newtheorem*{conjectureNN}{Conjecture}
\newcommand{\bconjNN}{\begin{conjectureNN}}
\newcommand{\econjNN}{\end{conjectureNN}}

\newtheorem{example}{Example}
\newcommand{\bex}{\begin{example}}
\newcommand{\eex}{\end{example}}

\newtheorem{problem}{Problem}
\newcommand{\bprob}{\begin{problem}}
\newcommand{\eprob}{\end{problem}}

\newtheorem*{problemNN}{Problem}
\newcommand{\bprobNN}{\begin{problemNN}}
\newcommand{\eprobNN}{\end{problemNN}}

\newtheorem{oproblem}{Open Problem}
\newcommand{\boprob}{\begin{oproblem}}
\newcommand{\eoprob}{\end{oproblem}}

\newtheorem*{oproblemNN}{Open Problem}
\newcommand{\boprobNN}{\begin{oproblemNN}}
\newcommand{\eoprobNN}{\end{oproblemNN}}

\newcommand{\beq}{\begin{equation}}
\newcommand{\eeq}{\end{equation}}
\newcommand{\benum}{\begin{enumerate}}
\newcommand{\eenum}{\end{enumerate}}

\newcommand{\Q}{\ensuremath{ \mathbf{Q} }}

\newcommand{\R}{\ensuremath{\mathbf R}}

\newcommand{\mbc}{\ensuremath{\mathbf c}}

\newcommand{\mbr}{\ensuremath{\mathbf r}}

\newcommand{\bq}{\begin{eqnarray*}}
\newcommand{\eq}{\end{eqnarray*}}
\newcommand{\be}{\begin{eqnarray}}
\newcommand{\ee}{\end{eqnarray}}
\newcommand{\ba}{\begin{array}}
\newcommand{\ea}{\end{array}}
\newcommand{\bfr}{\begin{flushright}}
\newcommand{\efr}{\end{flushright}}

\newcommand{\bmat}{\left(\begin{matrix}}
\newcommand{\emat}{\end{matrix}\right)}
\newcommand{\bsmallmat}{\left(\begin{smallmatrix}}
\newcommand{\esmallmat}{\end{smallmatrix}\right)}

\DeclareMathOperator{\col}{\text{col}}

\DeclareMathOperator{\diag}{\text{diag}}

\DeclareMathOperator{\qqand}{\qquad\text{and}\qquad}

\DeclareMathOperator{\row}{\text{row}}

\usepackage[all]{xy}
\usepackage{amssymb,latexsym}

\title[Alternate minimization]{Alternate minimization and doubly stochastic matrices}
\author{Melvyn B. Nathanson}\address{Department of Mathematics\\Lehman College (CUNY)\\Bronx, NY 10468}\email{melvyn.nathanson@lehman.cuny.edu}

\subjclass[2010]{11C20, 11B75,11B99.}

\keywords{Alternate minimization, Sinkhorn limits, diophantine approximation.}

\date{\today}

\begin{document}
\maketitle

\begin{abstract}
Sinkhorn's alternative minimization algorithm applied to a positive $n\times n$ matrix 
converges to a doubly stochastic matrix.  
If the algorithm, applied to a $2\times 2$ matrix, converges in a finite number of iterations, 
then it converges in at most two iterations, and the structure of such matrices is determined.  
\end{abstract}

\section{The alternate minimization algorithm}

A \emph{positive matrix}  is a matrix with positive coordinates.
Let $\diag(x_1,\ldots, x_n)$ denote the $n\times n$ diagonal matrix 
with coordinates $x_1,\ldots, x_n$ on the main diagonal.
A  \emph{positive diagonal matrix} is a diagonal matrix whose diagonal coordinates are positive.  
If $A$ is an $m\times n$  positive matrix,  
$X$ is an $m\times m$ positive diagonal matrix, 
and $Y$ is an $n\times n$ positive diagonal matrix, 
then $XA$ and $AY$ are $m\times n$ positive matrices. 

Let $A = (a_{i,j})$ be an $n \times n$  matrix.
The $i$th \emph{row sum} of $A$ is 
\[
\row_i(A) = \sum_{j=1}^n a_{i,j}. 
\]
The $j$th \emph{column sum} of $A$ is 
\[
\col_j(A)  = \sum_{i=1}^n a_{i,j}.
\]
The matrix $A$ is \emph{row stochastic} if $\row_i(A) = 1$ 
for all $i \in \{1,\ldots, n\}$.  
The matrix $A$ is \emph{column stochastic} if $\col_j(A) = 1$ 
for all $j \in \{1,\ldots, n\}$.  
The matrix $A$ is \emph{doubly stochastic} if it is 
both row stochastic and column stochastic.  

For example, a positive $2\times 2$  matrix $A$ is doubly stochastic 
if and only if there exist $\alpha, \beta  \in (0,1)$ such that  $\alpha + \beta = 1$ and 
\[
A = \bmat
\alpha &  \beta \\
 \beta & \alpha
\emat.
\]
If $\alpha, \beta, \gamma \in (0,1)$ satisfy $\alpha + 2\beta = \beta + 2\gamma = 1$, then 
the $3\times 3$ symmetric matrix 
\[
\bmat
\alpha & \beta & \beta \\
\beta & \gamma & \gamma \\
\beta & \gamma & \gamma
\emat 
\]
is doubly stochastic.

Let $A = (a_{i,j})$ be a positive $n \times n$  matrix.
We have $\row_i(A) > 0$ and $\col_j(A)>0$ for all $i,j \in \{1,\ldots, n\}$.  
Define the $n \times n$ positive diagonal matrix 
\[
X(A) = \diag \left( \frac{1}{\row_1(A)}, \frac{1}{\row_2(A)},\ldots, \frac{1}{\row_n(A)} \right).
\]
Multiplying $A$ on the left by $X(A)$ multiplies each coordinate in the $i$th 
row of $A$ by $ 1/\row_i(A)$, and so 
\begin{align*}
\row_i\left( X(A) A\right) 
& = \sum_{j=1}^n (X(A) A)_{i,j} 
 = \sum_{j=1}^n \frac{a_{i,j}}{\row_i(A)} 
 = \frac{ \row_i(A)}{\row_i(A)}  = 1
\end{align*}
for all $i \in \{1,2,\ldots, n\}$.  
The process of multiplying $A$ on the left by $X(A)$ to obtain the 
row stochastic matrix $X(A) A$ is called 
\emph{row scaling} or \emph{row normalization}.  
We have $X(A) A = A$   if and only if $A$  is row stochastic if and only if  $X(A) = I$.  
Note that the row stochastic matrix $X(A)A$ is not necessarily column stochastic.

Similarly, we define the $n \times n$ positive diagonal matrix 
\[
Y(A) = \diag \left( \frac{1}{\col_1(A)}, \frac{1}{\col_2(A)},\ldots, 
\frac{1}{\col_n(A)} \right).
\]
Multiplying $A$ on the right by $Y(A)$ multiplies 
each coordinate in the $j$th column 
of $A$ by $1/\col_j(A)$, and so 
\[
\col_j(AY(A)) = \sum_{i=1}^n (A Y(A))_{i,j} = \sum_{i=1}^n  \frac{a_{i,j}}{\col_j(A)} 
= \frac{\col_j(A)}{\col_j(A)} = 1
\]
for all $j \in \{1,2,\ldots, n\}$.  
The process of multiplying $A$ on the right by $Y(A)$ to obtain a 
column stochastic matrix $A Y(A)$ is called 
\emph{column scaling} or \emph{column normalization}.  
We have $AY(A) = A$ if and only if  $Y(A) = I$ if and only if $A$ is column stochastic.  
The column stochastic matrix 
$A Y(A)$ is not necessarily row stochastic.  

The following elementary identity shows that column scaling can be replaced 
by row scaling, and conversely.  .

\bl
Let $A^t$ denote the transpose of the $n \times n$ positive matrix $A = (a_{i,j})$.   
Row and column scaling satisfy the following transpose symmetries:
\[
AY(A) = \left( X(A^t) \ A^t \right)^t
\]
and
\[
X(A) A =  \left( A^t \ Y(A^t) \right)^t.
\]
\el

\begin{proof}
Let $A^t = (a_{i,j}^t)$, where $a_{i,j}^t = a_{j,i}$.  We have 
\[
\row_i(A^t)  = \sum_{j=1}^n a_{i,j}^t  = \sum_{j=1}^n a_{j,i}  = \col_i(A) 
\]
and so 
\begin{align*}
X(A^t) & = \diag\left( \frac{1}{\row_1(A^t)}, \ldots, \frac{1}{\row_n(A^t)} \right) \\
& = \diag\left( \frac{1}{\col_1(A)}, \ldots, \frac{1}{ \col_n(A)} \right) \\ 
& = Y(A).  
\end{align*}
Because the transpose of a diagonal matrix $D$ is $D^t = D$, we obtain 
\[
 \left( X(A^t) \ A^t \right)^t 
=  \left( \ A^t \right)^t  \left( X(A^t)\right)^t  
 = A \  X(A^t) =  A \  Y(A).
\]
The proof of the identity $X(A) A =  \left( A^t \ Y(A^t) \right)^t$ is similar.
\end{proof}

For example, if $A = \bmat a & b \\ c & d \emat$, 
then $A^t  =  \bmat a & c \\ b & d \emat$ and 
\[
X(A^t)  = \bmat 1/(a+c) & 0 \\ 0 & 1/(b+d) \emat = Y(A).
\]
We have 
\begin{align*}
 \left( X(A^t) \ A^t \right)^t 
 & =  \bmat a/(a+c) & c/(a+c)  \\ b/(b+d) & d/(b+d) \emat^t \\
 & =  \bmat a/(a+c) &  b/(b+d) \\ c/(a+c)  & d/(b+d) \emat  \\
 & = A \ Y(A).
\end{align*}

Sinkhorn~\cite{sink64} proved that row and column scaling satisfy the following uniqueness theorem. 

\bt           \label{SinkhornFinite:theorem:SinkhornUnique}
Let $A$ be a positive matrix, and let $X_1$, $X_2$, $Y_1$, and $Y_2$ 
be positive diagonal matrices.  
If 
\[
S_1 = X_1AY_1 \qqand S_2 = X_2 AY_2
\]
are doubly stochastic matrices, then 
\[
S_1 = S_2
\]
and there exists $\lambda > 0$ 
such that 
\[
X_2 = \lambda X_1 \qqand Y_2 = \lambda^{-1} Y_1.
\]
If the positive matrix $A$ is symmetric, then there is a unique positive diagonal matrix $D$ 
such that $S = DAD$ is doubly stochastic.  
\et

The following algorithm is called  ``alternate minimization'' 
(perhaps, more appropriately called ``alternate scaling'' 
or ``alternate normalization'').  
The proof of the convergence of the algorithm is due to  
Sinkhorn~\cite{sink64} and Sinkhorn and Knopp~\cite{sink-knop67}.

\bt     \label{SinkhornFinite:theorem:Sinkhorn}
Let $A = (a_{i,j})$ be a positive $n \times n$ matrix.
Construct inductively an infinite sequence of positive $n \times n$ 
matrices by alternate operations of column scaling and row scaling: 
\begin{align*}
A^{(0)}& = A \\
A^{(1)} & = A^{(0)}  \ Y\left( A^{(0)} \right)  \\ 
A^{(2)} & = X\left( A^{(1)} \right)   \  A^{(1)} \\
A^{(3)} & = A^{(2)}  \  Y\left( A^{(2)} \right)  \\ 
A^{(4)} & = X\left( A^{(3)} \right) \  A^{(3)} \\
A^{(5)} & = A^{(4)}  \ Y\left( A^{(4)} \right)  \\ 
& \vdots
\end{align*}
The sequence of matrices $\left( A^{(\ell)} \right)_{\ell=0}^{\infty}$ 
converges to a doubly stochastic matrix $S(A)$, 
and there exist positive diagonal matrices $X$ and $Y$ such that 
\[
S(A) = X A  Y.  
\]
\et

The sequence of matrices  $\left( A^{(\ell)} \right)_{\ell=0}^{\infty}$ is called the 
\emph{alternate minimization sequence}\index{alternate minimization sequence}  associated with $A$, 
and the matrix 
\[
S(A) = \lim_{\ell\rightarrow \infty} A^{(\ell)}
\]
is the \emph{alternate minimization limit}\index{alternate minimization limit} 
(also called the  \emph{Sinkhorn limit}\index{Sinkhorn limit}) of $A$.  

For example, if 
\[
A = A^{(0)} = \bmat 1 & 3 \\ 3 & 4 \emat
\]
then the next three matrices in the sequence  $\left( A^{(\ell)} \right)_{\ell=0}^{\infty}$ are 
\begin{align*}
 A^{(1)} & =  A^{(0)}  \ Y\left( A^{(0)} \right) = \bmat 1/4 & 3/7 \\ 3/4 & 4/7 \emat \\ 
 A^{(2)} & = X\left( A^{(1)} \right)   \  A^{(1)} = \bmat 7/19 & 12/19 \\  21/37 & 16/37  \emat \\
 A^{(3)} & = A^{(2)}  \  Y\left( A^{(2)} \right) =  \bmat 37/94 & 111/187 \\ 57/94 & 76/187 \emat 
\end{align*}

Let $P$ and $Q$ be positive diagonal $n \times n$ matrices.  
It follows from Theorem~\ref{SinkhornFinite:theorem:SinkhornUnique} 
that the alternate minimization limit of the positive $n \times n$ matrix $A$ 
is equal to the alternate minimization limit of the  matrix $PAQ$.
In particular, the matrices $A$ and $X(A) A$ have the same limits, 
and so  it makes no difference if we start the alternate minimization sequence 
by column scaling or by row scaling.

Let $A$ be a positive matrix, and let  $\left( A^{(\ell)} \right)_{\ell=0}^{\infty}$ 
be the alternate minimization sequence of matrices constructed 
in Theorem~\ref{SinkhornFinite:theorem:Sinkhorn}.  
If $A^{(L)}$ is doubly stochastic for some $L$, 
then $A^{(\ell)} = A^{(L)}$ for all $\ell \geq L$, 
and so the sequence of matrices $\left( A^{(\ell)} \right)_{\ell=0}^{\infty}$  
is eventually constant. In this presumably exceptional case, 
we say that the alternate minimization algorithm 
terminates in at most $L$ steps.  
Note that, if the $n \times n$ matrix $A$ has positive rational coordinates, 
then the matrix $A^{(\ell)}$ has positive rational coordinates for all $\ell \geq 1$.  
It follows that, if the Sinkhorn limit has irrational coordinates, 
then the alternate minimization algorithm cannot terminate in a finite number of steps.  
In Section~\ref{SinkhornFinite:section:finite2x2}, we prove that, 
for $2\times 2$ matrices, if the algorithm terminates in a  finite number of steps, 
then the algorithm terminates in at most two steps.

There is a vast literature on alternate minimization algorithms and Sinkhorn limits.  
For a recent survey, see Idel~\cite{idel16}.  In complexity theory, it is the asymptotics of the 
approximating sequence $\left( A^{(\ell)} \right)_{\ell=0}^{\infty}$  
that is important (for example, Allen-Zhu,  Li,  Oliveira,  and Wigderson~\cite{alle-li-oliv-wigd17}).    
This paper is concerned with number theoretic aspects 
of the algorithm, and with the classification of matrices for which the 
alternate minimization algorithm terminates in a finite number of steps.  
It is also of interest to consider the application of the algorithm to simultaneous approximation 
of irrational numbers by rational numbers.

\section{Alternate minimization limits for $2 \times 2$ matrices}
\label{SinkhornFinite:section:limit2x2}

\bt                  \label{SinkhornFinite:theorem:2x2}
Let 
\[
A = \bmat a & b \\ c & d \emat
\]
be a positive $2\times 2$ matrix.
Define the positive diagonal matrices  
\[
X = \bmat \sqrt{cd} & 0 \\ 0 & \sqrt{ab} \emat
\]
and
\[
Y = \bmat 
\left(a\sqrt{cd} + c \sqrt{ab} \right)^{-1}  & 0 \\ 0 & \left(b\sqrt{cd} + d\sqrt{ab} \right)^{-1} 
\emat
\]
The limit of the alternate minimization sequence $\left( A^{(\ell)} \right)_{\ell=0}^{\infty}$ 
is the doubly stochastic matrix 
\beq                             \label{SinkhornFinite:2x2A}
S(A) = XAY = \bmat \alpha & \beta  \\ \beta  & \alpha \emat 
\eeq
with 
\beq                             \label{SinkhornFinite:2x2B}
\alpha = \frac{\sqrt{ad}}{\sqrt{ad}+\sqrt{bc}}
\qqand 
\beta  = \frac{\sqrt{bc}}{\sqrt{ad}+\sqrt{bc}}.  
\eeq
\et

\begin{proof}
Simply compute the product $XAY$.  
That the matrix $XAY$ is the alternate minimization  
limit follows from  uniqueness (Theorem~\ref{SinkhornFinite:theorem:SinkhornUnique}).  
\end{proof}

\bc
Let $A = \bmat a & b \\ c & d \emat \in M_2^+(\Q)$.  
The alternate minimization limit of the matrix $A$ 
has rational coordinates if and only if $ad/bc$ is the square of a rational number.  
\ec

For example, if $A_1 =  \bmat 1 & 3 \\ 3 & 4 \emat $, then  
\begin{align*}
S(A_1) 
& = \bmat 2\sqrt{3} & 0 \\ 0 & \sqrt{3} \emat 
\bmat 1 & 3 \\ 3 & 4\emat 
\bmat \left( 5\sqrt{3} \right)^{-1} & 0 \\ 0 & \left( 10\sqrt{3} \right)^{-1} \emat \\ 
& = \bmat 2 & 0 \\ 0 & 1\emat 
\bmat 1 & 3 \\ 3 & 4 \emat 
\bmat  1/5 & 0 \\ 0 & 1/10 \emat \\ 
& = \bmat  2/5 & 3/5 \\ 3/5 & 2/5  \emat.
\end{align*} 

If $A_2 = \bmat 1 & 2 \\ 3 & 4 \emat $, then  
\begin{align*}
S(A_2) 
& =  \bmat 2\sqrt{3} & 0 \\ 0 & \sqrt{2} \emat 
\bmat 1 & 2 \\ 3 & 4 \emat 
\bmat \left(2\sqrt{3} + 3 \sqrt{2} \right)^{-1}  & 0 \\ 0 & \left(4\sqrt{3} + 4\sqrt{2} \right)^{-1}  \emat \\ 
& = \bmat \sqrt{6}-2 & 3 - \sqrt{6} \\ 3 - \sqrt{6}  & \sqrt{6}-2  \emat .
\end{align*} 
Because $A_2$ has rational coefficients and $S(A_2)$ has irrational coefficients, 
the alternate minimization algorithm for $A_2$ must have infinite length,  
that is, does not terminate in a finite number of steps.

\bt                         \label{SinkhornFinite:theorem:2x2symmetric}
Consider the  positive symmetric matrix 
\[
A = \bmat a & b \\ b & d \emat.
\]
Let 
\[
\lambda = \left( abd + b^2 \sqrt{ad} \right)^{-1/2}
\]
and 
\[
D =  \bmat \lambda \sqrt{bd} & 0 \\ 0 & \lambda \sqrt{ab} \emat.
\]
The Sinkhorn limit  of $A$ is the doubly stochastic matrix 
\[
S(A) = DAD = \bmat \alpha & \beta  \\ \beta  & \alpha \emat
\]
with 
\[
\alpha = \frac{\sqrt{ad}}{\sqrt{ad}+ b } 
\qqand 
\beta = \frac{b}{\sqrt{ad}+ b }.
\]
\et

\begin{proof}
The row scaling matrix  
\[
X(A) = \bmat \sqrt{bd} & 0 \\ 0 & \sqrt{ab} \emat
\]
and. the column scaling 
\[
Y(A) = \bmat 
\left( a\sqrt{bd} + b \sqrt{ab} \right)^{-1}  & 0 \\ 0 & \left( b\sqrt{bd} + d\sqrt{ab} \right)^{-1} 
\emat
\]
satisfy 
\[
D = \lambda X(A) = \lambda^{-1} Y(A).
\]
By Theorem~\ref{SinkhornFinite:theorem:2x2}, the matrix 
\[
XAY = (\lambda X)A \left( \lambda^{-1}Y \right) 
= DAD = \bmat \alpha & \beta \\ \beta & \alpha \emat 
\]
is doubly stochastic with  $\alpha = \sqrt{ad}/(\sqrt{ad}+b)$.
This completes the proof.  
\end{proof}

For example, if $A_3 = \bmat 1 & 2 \\ 2 & 4 \emat $, then  
\[
D = \bmat \sqrt{2}/2 & 0 \\ 0 & \sqrt{2}/4 \emat
\]
and
\begin{align*}
DA_3D  
& = \bmat \sqrt{2}/2 & 0 \\ 0 & \sqrt{2}/4 \emat
\bmat 1 & 2 \\ 2 & 4 \emat 
\bmat \sqrt{2}/2 & 0 \\ 0 & \sqrt{2}/4 \emat \\ 
& = \bmat  1/2 & 1/2 \\ 1/2 & 1/2  \emat.
\end{align*} 
If $A_4 = \bmat 1 & 1 \\ 1 & 1 \emat $, then  
\[
D = \bmat  1/\sqrt{2} & 0 \\ 0 &1/\sqrt{2} \emat
\]
and
\begin{align*}
DA_4D  
& = \bmat 1/\sqrt{2}  & 0 \\ 0 &1/\sqrt{2}  \emat
\bmat 1 & 1 \\ 1 & 1 \emat 
\bmat 1/\sqrt{2}  & 0 \\ 0 & 1/\sqrt{2} \emat \\ 
& = \bmat  1/2 & 1/2 \\ 1/2 & 1/2  \emat.
\end{align*} 
Note that the matrices $\bmat 1 & 2 \\ 2 & 4 \emat$ and 
$\bmat 1 & 1 \\ 1 & 1 \emat$ have the same Sinkhorn limits.


\section{Limits for $2 \times 2$ matrices in finitely many steps}
\label{SinkhornFinite:section:finite2x2}

\bt                         \label{SinkhornFinite:theorem:OneStep}
Let $A$ be a positive $2\times 2$ matrix that is not doubly stochastic.  
If the column scaled matrix $A Y(A)$ is doubly stochastic, then $A$ is a matrix of the form 
\beq                      \label{SinkhornFinite:OneStep-column}
A = \bmat a & ct \\ c & at \emat
\eeq
and
\[
S(A) = A Y(A) = \bmat a/(a+c) & c/(a+c) \\ c/(a+c) & a/(a+c) \emat.
\]
If the row scaled matrix $X(A) A$ is doubly stochastic, then $A$ is a matrix of the form 
\beq                      \label{SinkhornFinite:OneStep-row}
A = \bmat a & b \\ bt & at \emat
\eeq
and
\[
S(A) = X(A) A = \bmat a/(a+b) & b/(a+b) \\ b/(a+b) & a/(a+b) \emat.
\]
\et

For example, column scaling the matrix 
$\bmat 1 & 12 \\ 3 & 4  \emat$
and row scaling the matrix 
$\bmat 1 & 3 \\ 12 & 4  \emat$
both produce the doubly stochastic matrix $\bmat 1/4 & 3/4 \\ 3/4 & 1/4  \emat$.

\begin{proof}
Column scaling a matrix of the form~\eqref{SinkhornFinite:OneStep-column}
and row scaling a matrix of the form~\eqref{SinkhornFinite:OneStep-row} 
both produce doubly stochastic matrices.

Conversely, let $A = \bmat a & b \\ c & d \emat$.  
The column scaled matrix 
\[
AY(A) = \bmat a/(a+c) & b/(b+d) \\ c/(a+c) & d/(b+d) \emat
\]
is doubly stochastic if and only if 
\[
\frac{a}{a+c} + \frac{b}{b+d} = \frac{c}{a+c} + \frac{d}{b+d} = 1
\]
if and only if 
\[
ab = cd.
\]
Defining $t = b/c = d/a$, we obtain 
\[
A = \bmat a & ct \\ c & at \emat 
\qqand
S(A) = AY(A) = \bmat a/(a+c) & c/(a+c)  \\
c/(a+c) & a/(a+c) 
  \emat.
\]

Similarly, the row scaled matrix 
\[
 X(A) A = \bmat a/(a+b) & b/(a+b) \\ c/(c+d) & d/(c+d) \emat
\]
is doubly stochastic if and only if 
\[
\frac{a}{a+b} + \frac{c}{c+d} =  \frac{b}{a+b} + \frac{d}{c+d} = 1
\]
if and only if 
\[
ac = bd.
\]
Defining $t = c/b = d/a $, we obtain 
\[
A = \bmat a & b \\ bt & at \emat 
\qqand
S(A) = X(A) A = \bmat
a/(a+b) & b/(a+b) \\
b/(a+b) & a/(a+b)
\emat.
\]
This completes the proof.  
\end{proof}

\bt      \label{SinkhornFinite:theorem:FiniteAlgorithm}
Let $A$ be a positive $2\times 2$ row stochastic matrix that is not column stochastic. 
If column scaling $A$ produces a doubly stochastic matrix $S(A) = AY(A)$, 
then $A = \bmat a & 1-a \\ a & 1-a \emat$ with $a\neq 1/2$, 
and $S(A) = \bmat 1/2 & 1/2 \\ 1/2 & 1/2 \emat$.

Let $A$ be a positive $2\times 2$ column stochastic matrix that is not row stochastic.
If row scaling $A$ produces a doubly stochastic matrix $S(A) = X(A) A$, 
then $A = \bmat a & a \\ 1-a & 1-a \emat$ with $a\neq 1/2$, and 
$S(A) = \bmat 1/2 & 1/2 \\ 1/2 & 1/2 \emat$.
\et

\begin{proof}
Let $A$ be a positive $2\times 2$ matrix that is not doubly stochastic. 
By Theorem~\ref{SinkhornFinite:theorem:OneStep}, 
if column scaling $A$ produces a doubly stochastic matrix, 
then $A = \bmat a & ct \\ c & at \emat$  for some $t > 0$. 
If $A$ is also row stochastic, then 
\[
a+ct = c+at = 1
\]
and so 
\[
 (a-c)(1-t)=0.
\]
If $t=1$, then $a+c=1$ and $A$ is doubly stochastic, which is absurd.
Therefore, $t \neq 1$ and $a=c$.   It follows that 
$A = \bmat a & at \\ a & at \emat = \bmat a & 1-a \\ a & 1-a \emat$ with $a\neq 1/2$  
and $AY(A) = \bmat 1/2 & 1/2 \\ 1/2 & 1/2 \emat$.  
 
Similarly, if row scaling $A$ produces a doubly stochastic matrix, 
then Theorem~\ref{SinkhornFinite:theorem:OneStep} implies that 
$A = \bmat a & b \\ bt & at \emat$  for some $t > 0$. 
If $A$ is also column stochastic, then 
\[
a+bt = b+at = 1
\]
and so 
\[
 (a-b)(1-t)=0.
\]
If $t=1$, then $a+b=1$ and $A$ is doubly stochastic, which is absurd.
Therefore, $t \neq 1$ and $a=b$.   It follows that  
$A = \bmat a & a \\ at & at \emat = \bmat a & a \\ 1-a & 1-a \emat$ with $a\neq 1/2$  
and $X(A)A = \bmat 1/2 & 1/2 \\ 1/2 & 1/2 \emat$.  
This completes the proof.  
\end{proof}

\bt                    \label{SinkhornFinite:theorem:FiniteAlgorithm-2}
Let $A$ be a positive $2\times 2$ matrix that is not doubly stochastic.  
If the alternate minimization algorithm produces 
a doubly stochastic matrix $S(A)$ in a finite number of steps,   
then the algorithm terminates in at most two steps. 

Suppose that the algorithm terminates in exactly two steps.
The matrix $S(A)  = A^{(2)}$ is obtained from $A^{(1)}$ by column scaling if and only if 
there exist positive real numbers $p$, $r$, and $t$ with $t \neq 1$ such that 
\[
A = \bmat p & pt \\ r & rt \emat. 
\]
The matrix $S(A)  = A^{(2)}$ is obtained from $A^{(1)}$ by row scaling if and only if 
there exist positive real numbers $p$, $q$, and $t$ with $t \neq 1$ such that 
\[
A = \bmat p & q \\ pt & qt \emat. 
\]
In both cases, the alternate minimization limit is 
$S(A) = \bmat  1/2 & 1/2 \\ 1/2 & 1/2 \emat$.
\et

Note that we obtain the limit matrix $S(A)$ either by first column scaling and then row scaling, 
or by first row scaling and then column scaling.

\begin{proof}  
Let $L$ be a positive integer such that the alternate minimization algorithm for $A$ 
terminates in exactly $L$ steps.   There is a sequence of matrices $\left( A^{(\ell}) \right)_{\ell = 0}^L$ 
with $A^{(0)} = A$ and $A^{(L)} = S(A)$ such that, for $\ell = 1,\ldots, L$, 
the matrix  $A^{(\ell}$ is obtained from $A^{(\ell-1}$ 
by alternate column and row scalings.

Suppose that $L \geq 3$.  There are two cases.  
Either $A^{(L)}$ is obtained from $A^{(L-1)}$ by column scaling,  
or $A^{(L)}$ is obtained from $A^{(L-1)}$ by row scaling,  

If $A^{(L)}$ is obtained from $A^{(L-1)}$ by column scaling,  
then $A^{(L-1)}$ is obtained from $A^{(L-2)}$ by row scaling, 
and $A^{(L-2)}$ is obtained from $A^{(L-3)}$ by column scaling.  
We have the diagram
\[
\xymatrix{
A = A^{(0)} \ar[r] & ...\ar[r]   & A^{(L-3)} \ar[r]^{\text{col}}  & A^{(L-2)} \ar[r]^{\text{row}}  & A^{(L-1)}\ar[r]^{\text{col}}  & A^{(L)} = S(A).
}
\]
The matrix $A^{(L-1)} = X(A^{(L-2)} ) A^{(L-2)}$ is row stochastic but not column stochastic.  
By Theorem~\ref{SinkhornFinite:theorem:FiniteAlgorithm}, 
$A^{(L-1)} = \bmat a & 1-a \\ a & 1-a \emat$ with $a\neq 1/2$.  
If $L \geq 3$, then the matrix $A^{(L-2)}$ is column stochastic, and 
$A^{(L-1)} = X(A^{(L-2)})A^{(L-2)}$.  We have
\[
A^{(L-2)} = \bmat u & v \\ 1-u & 1-v \emat
\]
for some $u,v \in (0,1)$, and 
\begin{align*}
\bmat a & 1-a \\ a & 1-a \emat 
& = A^{(L-1)} 
 = X(A^{(L-2)})A^{(L-2)} \\ 
& =  \bmat u/(u+v) & v/(u+v)  \\ (1-u)/(2-u-v)  & (1-v)/(2-u-v)  \emat.
\end{align*}
Therefore,
\[
\frac{u}{u+v} = a = \frac{1-u}{2-u-v}.
\]
Equivalently,
\[
2u -u^2-uv = u(2-u-v) = (1-u)(u+v) = u+v-u^2-uv
\]
and so $u=v$ and 
\[
A^{(L-2)} = \bmat u & u \\ 1-u & 1-u \emat.
\]
Thus, the matrix 
\begin{align*}
A^{(L-1)}  & = X(A^{(L-2)})A^{(L-2)} = \bmat 1/(2u) & 0 \\ 0 & 1/(2-2u) \emat  \bmat u & u \\ 1-u & 1-u \emat \\
& =  \bmat 1/2 &1/2  \\ 1/2  & 1/2   \emat 
\end{align*}
 is doubly stochastic, which is  absurd.  
 Therefore, $A^{(L-2)}$ is  not column stochastic, and so $L \leq 2$.  
 
Suppose that $L=2$.  
 Let $A = A^{(0)} = \bmat p & q \\ r & s \emat$.  
Because $A^{(1)}$ is row stochastic but not column stochastic and $A^{(2)}$ is doubly stochastic, 
there exists $a \in (0,1)$,  $a \neq 1/2$, such that  
\begin{align*}
\bmat a & 1-a \\ a & 1-a \emat 
& = A^{(1)} 
 = X\left(A^{(0)} \right) A^{(0)} \\ 
& =  \bmat p/(p+q) & q/(p+q)  \\ r/(r+s) & s/(r+s)  \emat
\end{align*}
and so 
\[
\frac{p}{p+q} = \frac{r }{r+s}.
\]
Equivalently, $ps = qr$ and $s = qr/p$.
Thus, with $t = q/p$, we obtain 
\[
A = A^{(0)} = \bmat p & q \\ r & qr/p \emat = \bmat p & pt \\ r & rt \emat.
\]
If $t=1$, then $A^{(1)}$ is doubly stochastic, which is absurd.  Therefore, $t \neq 1$.
Thus, if $L=2$, then the alternate minimization sequence is 
\[
A  = \bmat p & pt \\ r & rt \emat \rightarrow    \bmat 1/(1+t) & t/(1+t)   \\ 1/(1+t)  &t/(1+t)   \emat 
\rightarrow    \bmat 1/2 &1/2  \\ 1/2  & 1/2   \emat .
\]

A similar argument works in the second case, 
where the matrix $A^{(L)}$ is obtained from $A^{(L-1)}$ by row scaling.  
This completes the proof.  
\end{proof}


\section{An alternate minimization limit for an $n \times n$ matrix}

There are no  formulae analogous to~\eqref{SinkhornFinite:2x2A} and~\eqref{SinkhornFinite:2x2B} 
for the alternate minimization limit of a positive $3\times 3$  matrix.
Nathanson~\cite{nath19x} has explicitly computed the alternate minimization limits 
of some classes of symmetric positive $3\times 3$   matrices.  Here is a simple example 
of an explicit calculation.

Let $n \geq 3$ and $K > 0$.  
We consider the positive symmetric $n \times n$ matrix 
\[
A = \bmat
K & 1 & 1 & \cdots & 1 \\
1 & 1 & 1& \cdots & 1  \\
1 & 1 & 1 & \cdots & 1 \\
\vdots & & & &  \vdots   \\
1 & 1 & 1 & \cdots & 1 
\emat.
\]
By Theorems~\ref{SinkhornFinite:theorem:SinkhornUnique} and~\ref{SinkhornFinite:theorem:Sinkhorn}, 
there exists a unique positive diagonal matrix 
\[
D = \diag(x_1,x_2, \ldots, x_n) =  \bmat
x_1 & 0 & 0 & 0 & 0 \\
0 & x_2 & 0 & 0 & 0 \\
0 & 0 & x_3 & 0 & 0  \\
\vdots & & & &  \vdots   \\
0 & 0 & 0 & 0 & x_n
\emat
\]
such that the matrix 
\[
S(A) = DAD =  \bmat
Kx_1^2 & x_1 x_2 & x_1 x_3 & \cdots & x_1x_n \\
x_2 x_1 & x_2^2 & x_2 x_3 & \cdots & x_2 x_n\\
x_3 x_1 & x_3 x_2  & x_3^2  & \cdots & x_3x_n \\
\vdots & & & &  \vdots   \\
x_n x_1 & x_n x_2  & x_n x_3  & \cdots & x_n^2 
\emat 
\]
is doubly stochastic.  Equivalently, 
\[
Kx_1^2 + x_1 \sum_{j=2}^n x_j  = 1 
\]
and
\[
 x_i \sum_{j=1}^nx_j   = 1 
\]
for $i = 2,3,\ldots,n$.  
It follows that 
\[
x_i = \frac{1}{\sum_{j=1}^nx_j   }
\]
for $i = 2,3,\ldots,n$, and 
\[
S(A) = \bmat \alpha & \beta & \beta & \cdots & \beta \\ 
\beta & \gamma & \gamma & \cdots & \gamma \\
\beta & \gamma & \gamma & \cdots & \gamma \\
\vdots &&&& \vdots \\
\beta & \gamma & \gamma & \cdots & \gamma
\emat
\]
where 
\begin{align*}
\alpha & = K x_1^2 \\
\beta  & = x_1x_2 = \frac{1-\alpha}{n-1} \\ 
\gamma & = x_2^2 = \frac{1-\beta }{n-1} = \frac{n-2+ \alpha}{(n-1)^2}. 
\end{align*}
We obtain 
\[
\left(  \frac{1-\alpha}{n-1} \right)^2 = \beta^2 = x_1^2 x_2^2 = \frac{\alpha}{K} \left( \frac{n-2+ \alpha}{(n-1)^2} \right)
\]
and so 
\[
(K-1)\alpha^2 - (2K+n-2)\alpha + K = 0.
\]

If $K = 1$, then $\alpha = \beta = \gamma =  1/n$.

If $K \neq 1$, then 
\[
\alpha = \frac{2(K-1)+n \pm \sqrt{4(n-1)K+(n-2)^2}}{2(K-1)}.
\]
The inequality $0 < \alpha < 1$ implies that 
\[
\alpha = \frac{2(K-1)+n - \sqrt{4(n-1)K+(n-2)^2}}{2(K-1)}
\]
if $K > 1$ and if $0 < K < 1$.  

For example, if $n = 3$, then 
\[
\alpha = \frac{2K+1 - \sqrt{8K+1}}{2(K-1)}.
\]
If $n=3$ and $K=2$, then 
\[
\alpha = \frac{5-\sqrt{17}}{2} , \qquad 
\beta = \frac{-3+\sqrt{17}}{4} , \qquad 
\gamma =  \frac{7- \sqrt{17}}{8},
\]
\[
x_1 = \sqrt{\frac{5-\sqrt{17}}{4}} 
\qqand
x_2 = \frac{-3+\sqrt{17}}{\sqrt{5-\sqrt{17}}},
\]
and 
\[
S(A) = DAD = \bmat 
 \frac{5-\sqrt{17}}{2} &  \frac{-3+\sqrt{17}}{4} &  \frac{-3+\sqrt{17}}{4} \\
  \frac{-3+\sqrt{17}}{4} &  \frac{7 - \sqrt{17}}{8} &  \frac{7 - \sqrt{17}}{8} \\
 \frac{-3+\sqrt{17}}{4} &  \frac{7 - \sqrt{17}}{8} &  \frac{7 - \sqrt{17}}{8}
 \emat.
\]
If $n=3$ and $K=3$, then 
\[
\alpha  = \frac{1}{2}, \qquad 
\beta  = \frac{1}{4}, \qquad 
\gamma  =  \frac{3}{8}.
\]
For $n=3$ and  integers $K \geq 2$, the doubly stochastic matrix $S(A)$ is rational 
if and only if $K$ is a triangular number, that is, a number of the form $K = (k^2+k)/2$
for some positive integer $k$.
In this case, we have 
\begin{align*}
\alpha & = \frac{k^2-k}{k^2+k-2} \\
\beta & = \frac{k -1}{k^2+k-2} \\
\gamma & =   \frac{k^2-1}{2(k^2+k-2)}.
\end{align*}

If $n = 4$, then 
\[
\alpha = \frac{K+1 - \sqrt{3K+1}}{K-1}.
\]
If $n=4$ and $K=2$, then 
\[
\alpha =3-\sqrt{7} , \qquad 
\beta = \frac{-2+\sqrt{7}}{3} , \qquad 
\gamma =  \frac{5- \sqrt{7}}{9},
\]
If $n=4$ and $K=5$, then 
\[
\alpha = 1/2, \qquad 
\beta = 1/6 , \qquad 
\gamma = 5/18.
\]

\section{Open problems}

\bprob
Does there exist a positive $3\times 3$ matrix that is row stochastic but not column stochastic, 
and becomes doubly stochastic after one column scaling?  
This is equivalent to asking if there is a positive $3\times 3$ matrix that, with respect to the alternate minimization algorithm, 
has finite length $L \geq 2$.
\eprob

\bprob
 Let $n \geq 3$. Does there exist an integer $L^{\ast}(n)$ such that, 
 if $A$ is a positive $n\times n$ matrix for which the alternate minimization algorithm 
terminates in a finite number of steps, then the alternate minimization algorithm 
terminates in  at most $L^{\ast}(n)$ steps? 
\eprob

\bprob
Let $K$ be a subfield of \R, and let $M_n^{+}(K)$ be the set of positive $n\times n$ 
matrices with coordinates in $K$.  
If $A\in M_n^+(K)$, then $A^{(\ell)} \in M_n^+(K)$ for all matrices in the 
alternate minimization sequence $\left( A^{(\ell)} \right)_{\ell=0}^{\infty}$.  
It follows that if $S(A) \notin M_n^+(K)$, then the alternate minimization 
algorithm for the matrix $A$ has infinite length.
Thus, if $A \in M_n^+(\Q)$ and if the doubly stochastic limit $S(A)$ contains 
an irrational coordinate, then the alternate minimization 
algorithm has infinite length.  
In this case, the coordinates in the matrices $\left( A^{(\ell)} \right)_{\ell=0}^{\infty}$
are sequences of rational numbers that simultaneously converge to the coordinates of $S(A)$.  
It is of interest to understand the rate of convergence.  

\eprob

\bprob
Let $\mbr = \bmat r_1 \\ \vdots \\ r_m \emat \in \R^m$ and $\mbc = \bmat c_1 \\ \vdots \\ c_n \emat \in \R^n$ 
be vectors with positive coordinates such that 
\[
\sum_{i=1}^m r_i = \sum_{j=1}^n c_j.
\]
Let $A = (a_{i,j})$ be an $m\times n$ matrix.  
The matrix is $A$ is $\mbr$-row stochastic 
if 
\[
\row_i(A) = \sum_{j=1}^n a_{i,j} = r_i
\]
for all $i \in \{1,\ldots, m\}$.  
The matrix is $A$ is $\mbc$-column stochastic 
if 
\[
\col_j(A) = \sum_{i=1}^m a_{i,j} 
\]
for all $j \in \{1,\ldots, n\}$.
The matrix is $A$ is $(\mbr,\mbc)$ stochastic 
if it is both $\mbr$-row stochastic  and $\mbc$-column stochastic.

Let $A$ be a positive $m\times n$  matrix, and let 
\[
X_{\mbr}(A) = \diag \left( \frac{r_1}{\row_1(A)}, \frac{r_2}{\row_2(A)},\ldots, \frac{r_m}{\row_m(A)} \right) 
\]
and
\[
Y_{\mbc}(A) = \diag \left( \frac{c_1}{\col_1(A)}, \frac{c_2}{\col_2(A)},\ldots, 
\frac{c_n}{\col_n(A)} \right).
\]
The matrix $X_{\mbr}(A) \ A$ is $\mbr$-row stochastic, and 
the matrix $A \ Y_{\mbc} (A)$ is $\mbr$-column stochastic.  
The analogous $(\mbr,\mbc)$-alternate minimization algorithm applied to a positive $m\times n$ 
matrix always converges to an $(\mbr,\mbc)$-stochastic matrix.  

 Let $m,n \geq 2$. Does there exist an integer $L^{\ast}(m,n)$ such that, 
 if $A$ is a positive $m\times n$ matrix for which the 
 $(\mbr,\mbc)$-alternate minimization algorithm 
terminates in a finite number of steps, 
then the $(\mbr,\mbc)$-alternate minimization algorithm 
terminates in  at most $L^{\ast}(m,n)$ steps? 
\eprob

\bprob
Does there exist a constant $C_n$ with the following property:  
If $A$ is an positive $n \times n$ matrix such that the alternate minimization algorithm, 
starting with row scaling, terminates in $N_1$ steps, and the alternate minimization algorithm, 
starting with column scaling, terminates in $N_2$ steps, then $|N_1 - N_2| < C_n$?
\eprob

\textbf{Note added in proof.} 
S.~B. Ekhad and D.~Zeilberger 
(\emph{Answers to some questions about explicit Sinkhorn limits posed by Mel Nathanson}, 
arXiv:1902.10783) solved Problem 1 by constructing a positive $3\times 3$ matrix that is row stochastic but not column stochastic, and becomes doubly stochastic after one column scaling.
M. B. Nathanson (\emph{Matrix scaling limits in finitely many iterations}, arXiv:1903.06778) 
generalized this construction to $n \times n$ matrices.  

Alex Cohen (unpublished) solved Problem 2 by proving that $L^{\ast}(n) = 2$ for all $n \geq 3$.  
This also solves Problem 5.  Extending Cohen's proof, Nathanson (unpublished) solved Problem 4 
by showing that $L^{\ast}(m,n) = 2$ for all $m,n \geq 2$.

\def\cprime{$'$} \def\cprime{$'$} \def\cprime{$'$}
\providecommand{\bysame}{\leavevmode\hbox to3em{\hrulefill}\thinspace}
\providecommand{\MR}{\relax\ifhmode\unskip\space\fi MR }
\providecommand{\MRhref}[2]{%
  \href{http://www.ams.org/mathscinet-getitem?mr=#1}{#2}
}
\providecommand{\href}[2]{#2}

\end{document}